\documentclass[a4paper, 11pt]{amsart} 
\usepackage{graphicx}
\usepackage{psfrag}
\usepackage{amscd}

\usepackage{amsmath, amscd, amsthm, amssymb, mathrsfs}

\DeclareMathOperator{\im}{im}

\DeclareMathOperator{\res}{res}

\DeclareMathOperator{\U}{U}

\newcommand{\R}{\mathbb R}
\newcommand{\C}{\mathbb C}
\newcommand{\Z}{\mathbb Z}

\newcommand{\F}{\mathbb F}

\newcommand{\diff}{\text{\rm d}}
\newcommand{\del}{\partial}

\renewcommand{\P}{\mathbb P}

\theoremstyle{plain}
	\newtheorem{theorem}{Theorem}[section]
	\newtheorem{proposition}[theorem]{Proposition}
	\newtheorem{lemma}[theorem]{Lemma}

\theoremstyle{definition}

\theoremstyle{plain}
	\newtheorem*{theorem*}{Theorem}
	\newtheorem*{proposition*}{Proposition}
	\newtheorem*{lemma*}{Lemma}
	\newtheorem*{corollary*}{Corollary}
	\newtheorem*{conjecture*}{Conjecture}
\theoremstyle{definition}
	\newtheorem*{definition*}{Definition}
	\newtheorem*{remark*}{Remark}
	\newtheorem*{remarks*}{Remarks}

\numberwithin{equation}{section}

\begin{document}

\title{Toric anti-self-dual 4-manifolds via complex geometry}
\author{Simon Donaldson}
\author{Joel Fine}

\begin{abstract}
Using the twistor correspondence, this article gives a one-to-one correspondence between germs of toric anti-self-dual conformal classes and certain holomorphic data determined by the induced action on twistor space. Recovering the metric from the holomorphic data leads to the classical problem of prescribing the \v{C}ech coboundary of 0-cochains on an elliptic curve covered by two annuli. 

The classes admitting K\"ahler representatives are described; each such class contains a circle of K\"ahler metrics. This gives new local examples of scalar flat K\"ahler surfaces and generalises work of Joyce \cite{joyce} who considered the case where the distribution orthogonal to the torus action is integrable.
\end{abstract}

\maketitle

\section{Introduction and overview}

The twistor correspondence between anti-self-dual 4-manifolds and certain complex three-folds makes it possible to use techniques of complex geometry to solve problems in Riemannian geometry. This paper exploits this approach to describe the local geometry of an anti-self-dual 4-manifold which admits (the germ of) a conformally Killing $2$-torus action. (The relevant parts of the twistor correspondence are briefly reviewed below in Section \ref{review of twistors}.)

Our main result (Theorem \ref{main theorem}) relates the classification problem to some very simple holomorphic data: pairs $(\tau,\phi)$ where
\begin{itemize}
\item 
$\tau$ is a holomorphic involution of a neighbourhood of $0$ in $\C$ with $\tau'(0)=-1$;
\item 
$\phi$ is a holomorphic $\C^{2}$-valued function on the same
neighbourhood which is odd with respect to $\tau$ (i.e., $\phi(\tau(z))=-\phi(z)$) and such that $\phi'(0), \overline{\phi'(0)}$ are $\C$-linearly independent
 in $\C^{2}$. 
 \end{itemize}
Given such a pair we construct a twistor space, and hence an anti-self-dual $4$-manifold, with a local torus action. Conversely, we show that any such $4$-manifold which satisfies a mild genericity condition arises in this way.

Perhaps it should be noted that the neighbourhood mentioned above appears in the construction as an open set in $\P^1$. The natural symmetry group is then that of M\"obius transformations fixing a point rather than the much larger group of arbitrary holomorphic coordinate transformations, under which all the involutions $\tau$ are equivalent. In fact, the construction involves a real structure which corresponds to a choice of antipodal map on $\P^1$; the symmetry group is then those M\"obius transformations which fix a point and preserve this antipodal map, in other words rotations about a fixed axis.
 
There is already a considerable body of work on anti-self-dual $4$-manifolds with torus symmetries. Notably, Joyce \cite{joyce} gave a complete classification of {\it surface-orthogonal} solutions. Among other things, Joyce showed that that these are the same as the solutions containing a  Kahler metric in the conformal class for which the orbits of the torus action are Lagrangian. Joyce's approach is differential-geometric and does not use  twistor methods explicitly. His main conclusion is that the problem can be reduced to the solution of a certain linear PDE, which can be viewed as that defining axially-symmetric
 harmonic functions in $\R^{3}$. 
 
 To compare our results with those of Joyce, and for other purposes (including global questions which we do not go into at all here) it is important to have more explicit formulae for the metrics: that is, to implement the twistor correspondence explicitly.  We take this up in Section \ref{metric formulae}. We find that Joyce's surface-orthogonal hypothesis corresponds to the case when $\tau$ is the standard involution $\tau(z)=-z$. In this case, we are able to express the metric explicity in terms of certain contour integrals. These arise, from our point of view, in the solution of a Riemann--Hilbert problem which we formulate in terms of Cech cohomology on an elliptic curve. In turn, this Riemann--Hilbert problem arises in the description of the lines in our twistor space. In this way, we essentially achieve a twistor derivation of Joyce's results, through a  contour integral representation  of the solutions to the linear PDE.
 
In the general case, we find that the essential difficulty in writing the metric explicitly is that of uniformising a Riemann surface of genus $0$ which is constructed by gluing together patches using certain maps defined by $\tau$. If this can be done, the metric is given by contour integral formulae, much as before. We do not expect that this uniformisation problem can be solved explicitly in general but we find a family of cases where it can be. This gives rise to a family of anti-self-dual $4$-manifolds (which, as far as we know, are new) where the metric can be expressed explicitly in terms of contour integrals. In another direction, we discover a family of conformally-Kahler solutions (again, as far as we know, new) in which the orbits are not Lagrangian (Theorem \ref{conformally Kahler}). It would be interesting to solve the uniformisation problem, and hence find the metric, more explictly in this case, but we have not been able to do so yet.   

When this paper was close to completion we learnt from conversations with Lionel Mason that he and David Calderbank had obtained very similar results, in a project going back several years \cite{calderbank-mason}. While the principal results seem to be equivalent, we hope that the points of view and emphases are sufficiently different for  our contribution to be worthwhile.

\subsubsection*{Acknowledgements} The authors are grateful to Kevin Costello, Nigel Hitchin and Richard Thomas for helpful discussions.

\section{Review of the twistor correspondence}
\label{review of twistors}

The Penrose twistor correspondence is a one-to-one correspondence between conformal classes of anti-self-dual 4-manifolds $M$ and twistor spaces $Z$ which are complex three-folds with certain properties. In the context of Riemannian geometry this theory was first developed in detail by Atiyah--Hitchin--Singer \cite{atiyah-hitchin-singer}. The relevant parts are described briefly here.

Starting on the complex side, a twistor space is a complex three-fold $Z$ containing an embedded rational curve $L \subset Z$ with normal bundle $\nu_L \cong \mathcal O(1)\oplus\mathcal O(1)$. Now $H^1(\nu_L) = 0$ and $H^0(\nu_L)$ is four dimensional, so, by the deformation theory of embedded submanifolds, the set of deformations of $L$ inside $Z$ is parametrised by a complex four-fold $M^\C$. 

This four-fold carries a natural holomorphic conformal structure: a holomorphically varying, complex valued, symmetric, non-degenerate, bilinear form on each each tangent space, defined up to scale. Such a structure is determined by its null cone; the null vectors in $T_{L'}M^\C \cong H^0(\nu_{L'})$ are precisely those sections with a zero. In other words, two infinitesimally close lines in $Z$ are null separated if they intersect. 

The Weyl curvature of a conformal structure is defined in this complex setting in the same way as in the more familiar Riemannian context. As in four-dimensional Riemannian geometry, the Weyl curvature splits into self-dual and anti-self-dual parts. The main theorem on the complex side of the twistor correspondence is that the conformal class on $Z$ described above has anti-self-dual Weyl curvature.

In order to produce a four-dimensional Riemannian manifold an extra piece of data is needed, namely there should be a free anti-holomorphic involution $\gamma$ on $Z$ which preserves the inital line $L$. The involution $\gamma$ also acts on $M^\C$ with a real four-dimensional fixed set $M \subset M^\C$ corresponding to the real lines in $Z$. It follows from the fact that $\gamma$ is free that the complex conformal structure on $M^\C$ restricts to a Riemannian conformal structure on $M$ which is still anti-self-dual.

This construction of an anti-self-dual 4-manifold from a complex three-fold can be inverted. Given an oriented 4-manifold with a conformal structure, let $Z$ denote the bundle whose fibre over $x$ consists of almost complex structures on $T_xM$ compatible with the orientation and orthogonal with respect to the conformal class. $Z$ carries a natural almost complex structure which is integrable if and only if the conformal structure is anti-self-dual. The fibres of $Z$ are embedded rational curves with normal bundles $\mathcal O(1) \oplus \mathcal O(1)$. The antipodal map on each fibre defines an antiholomorphic involution on $Z$ recovering the complex picture described above. 

Many statements about the Riemannian geometry of an anti-self-dual 4-manifold have holomorphic interpretations on its twistor space. Two examples of this phenomenon which will be used later are the following.

\begin{theorem}[Pontecorvo \cite{pontecorvo}]
Let $Z$ be the twistor space of an anti-self-dual 4-manifold $M$. Then a K\"ahler representative of the conformal class is defined by a holomorphic section of $K_Z^{-1/2}$, compatible with $\gamma$ and non-zero on each fibre of $Z$.
\end{theorem}

\begin{proof}[Sketch proof]
Pontecorvo's original proof relies on Hitchin's general theory of linear field equations on anti-self-dual 4-manifolds \cite{hitchin}. It is also possible to give a direct proof via complex geometry, as sketched here.

Let $s\in H^0(K^{-1/2})$ be as in the hypotheses. The polar divisor of the meromorphic 3-form $\theta = s^{-2}$ intersects each line in $Z$ in two points, each with multiplicity two. Define a 2-form $\omega$ on the space of lines $M^\C$  as follows. Let $u,v \in T_LM^\C = H^0(\nu_L)$. The meromorphic 1-form on $L$ given by $\theta(u\wedge v)|_L$ has two double poles at  $J, J'$ say. Set
$$
\omega(u,v)
=
2\pi i \res_J \theta(u \wedge v)|_L.
$$

It is $(1,1)$ with respect to the complex structure $J$ since it vanishes when both $u$ and $v$ have zeros at $J$. It is positive on the real lines $M\subset M^\C$ as $\omega(u, Ju)$ positive whenever $u$ is non-zero and tangent to $M$.

To see that $\omega$ is closed, pick for each $L\in M^\C$ a circular contour $C_L$ on $L$ separating the points $J, J'$. Let 
$
Y = \{(L, q) : q\in C_L\}
$.
There are projections $p \colon Y \to M^\C$, $q \colon Y \to Z$. Then $\omega = p_* q^* \theta$. Since $\diff \theta = 0$ and $\diff$ commutes with $p_*$ and $q^*$, it follows that $\diff \omega =0$.
\end{proof}

In a similar vein, there is a twistor characterisation of hypercomplex structures on an anti-self-dual four manifold.

\begin{theorem}[Boyer \cite{boyer}]
Let $Z$ be the twistor space of an anti-self-dual 4-manifold $M$. Then a hypercomplex structure in the conformal class is defined by a holomorphic projection $\pi \colon Z \to \P^1$, compatible with $\gamma$ and which is an isomorphism on each fibre of $Z$.
\end{theorem}

\section{Twistor spaces with a $\C^2$-action}
\label{classification}
 
The remainder of this article focuses on the following situation: $(M, g)$ is an anti-self-dual 4-manifold admitting two linearly-independent, commuting, conformally Killing fields $X_1$ and $X_2$; by a slight abuse of terminology, such an $M$ is called toric. The Killing fields induce holomorphic vector fields $\tilde X_1$ and $\tilde X_2$ on the twistor space $Z$ of $M$ generating (the germ of) a $\C^2$-action. The idea is to use this $\C^2$-action to determine the complex geometry of $Z$ near a twistor line and hence, by the twistor correspondence, classify the local geometry of toric, anti-self-dual, 4-manifolds.  

\subsection{Non-trivial isotropy}

The first issue to address is whether or not the $\C^2$-action on $Z$ is free. The failure of the action to be free is measured by the zero locus $\Sigma$ of $\tilde X_1 \wedge \tilde X_2 \in H^0(\Lambda^2 TZ)$. As this is a section of a rank three vector bundle over a 3-fold, one might hope that this has empty intersection with a generic twistor line. This is not always the case, however.

For example, consider $\R^4$ with the Euclidean metric and Killing fields $X_i$ given by unit translation in two orthogonal  directions. The condition $JX_1 = X_2$ (along with metric and orientation compatibility) determines a complex structure $J$ on $\R^4$ and $\pm J$ are the only compatible complex structures for which $X_1$ and $X_2$ span a complex line. Since the $X_i$ are covariant constant, their holomorphic and horizontal lifts to $Z$ agree. This means the zero locus $\Sigma$ is the divisor in twistor space corresponding to $\pm J$. As is explained below, this phenomenon is due to the fact that the Euclidean metric is hyperk\"ahler and the Killing fields are triholomorphic. 

Return now to the general case. The horizontal components of the $\tilde X_i$ are the horizontal lifts of the $X_i$; on a twistor line $L$, they are linearly dependent at precisely the two antipodal points $\pm J$ for which $X_1$ and $X_2$ span a complex line. These are the only points at which  $\tilde X_1 \wedge \tilde X_2$ could possibly vanish and at every other point of $L$ the $\tilde X_i$ span a complex plane transverse to $L$. So $L \cap \Sigma$ is either empty or the pair $\pm J$.

The zero locus is $\C^2$-invariant. In particular, it can't be zero dimensional since both the $\tilde X_i$ are non-zero. There are three possibilities for $\Sigma$:

\begin{enumerate}
\item
$\Sigma = \emptyset$ and the $\C^2$-action is free. This is the generic case and the one which will be dealt with in this article.
\item
$\dim \Sigma = 1$. In this case, for a generic twistor line $L$, $L \cap \Sigma = \emptyset$. Near a generic line the $\C^2$-action is free as above.
\item
$\dim \Sigma = 2$. The following proposition shows that this case corresponds to a hypercomplex $M$ with $X_1$ and $X_2$ triholomorphic.
\end{enumerate}

\begin{proposition}\label{hypercomplex}
$\dim \Sigma = 2$ if and only if $(M, g)$ admits a compatible hypercomplex structure for which $X_1$ and $X_2$ are triholomorphic.
\end{proposition}

\begin{proof}
Begin by assuming that $\dim \Sigma = 2$. The aim is to find a holomorphic map $Z \to \P^1$ for which the real twistor lines are sections. To do this choose coordinates $\{z_1, z_2, z_3\}$ near $J$ in which $\tilde X_1 = \del_1$ and $L=(z_1 = 0 = z_2)$. Write $\tilde X_2 = \sum a_i \del_i$. Then $\Sigma =\{a_2 = 0 = a_3\}$. Since $\dim \Sigma = 2$, $a_2$ and $a_3$ have a common factor $c$; put $a_2= b c$. 

The projection of $\tilde X_1 \wedge \tilde X_2$ to $\Lambda^2 \nu_L$ vanishes at $\pm J$. Since $\Lambda^2 \nu_L \cong \mathcal O(2)$, these zeros are simple. This says that $\del_3 a_2 \neq 0$ at $J$. Since $c$ vanishes at $J$, this implies that $b$ is non-zero at $J$. The normal projections of the $\tilde X_i$ span the space
$
\left\langle
\del_1,
bc \del_2
\right\rangle
=
\left\langle
\del_1,
b \del_2
\right\rangle
$
which is a plane even at $J$ since $b$ is non-zero there. Hence the two-planes spanned by the $\tilde X_i$ extend to a transverse distribution over all of $L$.

Doing this on twistor lines near to $L$ defines a holomorphic two-plane distribution on $Z$ which is everywhere transverse to the lines. The distribution is integrable; the leaves through points other than $\pm J$ are the $\C^2$-orbits, whilst the leaves through $\pm J$ are the two components of $\Sigma$. The leaves define an identification of each nearby line with $L$ and so a map $Z \to L$. Moreover, since the $\tilde X_i$ are everywhere tangent to the distribution, they are triholomorphic with respect to the corresponding hypercomplex structure.

For the converse, assume that there is a map $Z \to \P^1$, with $\tilde X_i$ tangent to its fibres. Then $\Sigma$ is a pair of fibres of this map and hence has dimension 2.
\end{proof}

In the hyper\emph{k\"ahler} case, the local geometry is completely understood via the Gibbons--Hawking ansatz \cite{gibbons-hawking}. One can go a long way in analysing the general hypercomplex case, but in the remainder this paper we will focus on the generic case when the $\C^2$-action is free near a twistor line.

\subsection{Producing the holomorphic data}

Assume now that the $\C^2$-action on $Z$ is free near a twistor line $L$. It is possible to extract certain holomorphic data from the complex geometry of $Z$ near the line $L$ which completely characterises $Z$ and hence the local geometry of $M$.

Recall that the $\C^2$-orbits are transverse to $L$ except at the two antipodal points $\pm J$ for which $X_1$ and $X_2$ span a complex line. In other words, the section of $\Lambda^2 \nu_L$ obtained from projecting $\tilde X_1 \wedge \tilde X_2$ vanishes at $\pm J$. Since $\Lambda^2 \nu_L \cong \mathcal O(2)$ these zeros are simple, so the $\C^2$-orbits through $\pm J$ are tangential to $L$ to first order. Hence the orbit through a point $J'$ near to $J$ meets $L$ in one other nearby point $\tau(J')$. This defines a local holomorphic involution $\tau$ of a neighbourhood of $J$ in $L$. 

Since $J'$ and $\tau(J')$ lie on the same $\C^2$-orbit there is a unique vector $\phi(J') \in \C^2$ such that $\tau(J') = \phi(J') \cdot J'$. This determines a $\C^2$-valued holomorphic function $\phi$ defined on the domain of $\tau$ and which is odd with respect to $\tau$, meaning that $\phi \circ \tau = -\phi$. There is a distinguished line in $\C^2$, namely the one which runs tangentially to $L$ at $J$. This line is spanned by $\diff\phi_J$ which, consequently, must be non-zero.

There is a similar picture near $-J$, but this is related to that above by the real structure and so contains no extra information.   Note, however, that there are now two distinguished lines in $\C^2$, one tangential at $J$, the other at $-J$, and that these lines must be distinct: if not, the derivative of the $\C^2$-action in that direction would give a section of $\nu_L$ which vanishes twice contradicting $\nu_L \cong \mathcal O(1) \oplus \mathcal O(1)$. This translates into the condition $\im \diff \phi_J\neq\overline{\im \diff \phi_J}$.

In summary then, this procedure gives two pieces of holomorphic data:

\begin{itemize}
\item
A holomorphic involution $\tau$ defined on a neighbourhood $V$ of a point $0 \in \P^1$, which fixes only $0$.
\item 
A holomorphic function $\phi \colon V \to \C^2$ satisfying $\phi \circ \tau = -\phi$ and with $\im \diff \phi_0 \neq \overline{\im \diff \phi_0}$.
\end{itemize}

\subsection{Producing the twistor space}

Given such $\tau$ and $\phi$ it is possible to build a twistor space and hence a toric, anti-self-dual, 4-manifold. This construction is described below.

\subsubsection{The geometry near a point of tangency}

The first step is to model the geometry of the twistor space near a point where the $\C^2$-action is tangent to a twistor line.

Let $\tau$ and $\phi$ be as above. Choose a local coordinate, centred at $0\in \P^1$, in which $\tau(z) = -z$, and $\tau$ is defined on the disc $D$. Define a $\Z_2$-action on $D \times \C^2$ by $(z,v) \mapsto (-z, v + \phi(z))$ and write $Y = (D \times \C^2) / \Z_2$. 

The map $q \colon (z,v) \mapsto (z^2, v + \frac{1}{2} \phi (z))$ is $\Z_2$-invariant and shows $Y$ is biholomorphic to $D \times \C^2$. The central disc $D \times \{0\}$ projects under the quotient map $q$ to the curve $(z^2, \frac{1}{2} \phi (z))$. Since $\phi'(0) \neq 0$, this is an embedded disc $D' \subset Y$.

$\C^2$ acts on $D \times \C^2$ by translation in the second factor and this commutes with the $\Z_2$-action. Hence the $\C^2$-action descends to $Y$. Under the identification $Y \cong D\times \C^2$ this is just the standard $\C^2$-action again. It is transverse to the disc $D'$ everywhere except at the origin where the line in $\C^2$ spanned by $\phi'(0)$ is tangential to $D'$.

\subsubsection{Twistor space with the germ of a $\C^2$-action} 

Two copies of this construction can be ``glued together'' to produce a complex $3$-fold containing a smooth rational curve with normal bundle $\mathcal O(1) \oplus \mathcal O(1)$. This is done as follows.

Initially, no attempt will be made to produce a twistor space with a real structure; the anti-holomorphic involution will be introduced later. Accordingly, the construction begins with double the starting data, i.e., two locally defined involutions of $\P^1$, $\tau$ and $\sigma$, fixing $0$ and $\infty$ respectively, and two corresponding odd holomorphic functions $\phi$ and $\psi$, with $\diff \phi_0 \neq 0 \neq \diff \psi_\infty$.

Choose local coordinates near $0$ and near $\infty$ in which each of the involutions is given by $z \mapsto -z$; suppose, moreover, that $\tau$ and $\sigma$ are defined on discs $D_1, D_2 \subset \P^1$  respectively with $D_1 \cap D_2 = \emptyset$. Let $\Omega$ be an open set slightly larger than $\P^1 \setminus (D_1 \cup D_2)$ so that $D_1, D_2, \Omega$ is an open cover of $\P^1$. 

By shrinking the $D_i$ if necessary, assume that $\phi$ and $\psi$ are continuous and non-zero, except at the origins, on the closures $\overline D_i$. Let $V_1 \subset D_1$ be a smaller neighbourhood of $0$ and $V_2\subset D_2$ a smaller neighbourhood of $\infty$. There exists $\varepsilon> 0$ such that for all $z \in D_1 \setminus V_1$, $|\phi(z)| > \varepsilon$ and for all $z \in D_2 \setminus V_2$,  $|\psi(z)| > \varepsilon$. Let $U = \{ v \in \C^2 : |v| < \varepsilon\}$ and  define
$$
Y_1 = D_1 \times U / \sim \quad \text{where}\ (z, v) \sim (-z, v + \phi(z)),
$$
$$
Y_2 = D_2 \times U / \sim \quad \text{where}\ (z, v) \sim (-z, v + \psi(z)).
$$
As before, these quotients are complex manifolds and, since $\diff\phi_0$ and $\diff\psi_\infty$ are non-zero, each $Y_i$ contains an embedded copy of $D_i$, denoted $D'_i$. The germ of a $\C^2$-action on $D_i \times U$ descends to give the germ of an action on $Y_i$ which is transverse to $D'_i$ except at one point. In $Y_1$ the line $\im \diff \phi_0 \subset\C^2$ is tangential; in $Y_2$ the line $\im \diff\psi_\infty$ is tangential.

By the choice of $\varepsilon$, the quotient map $D_i \times U \to Y_i$ is an isomorphism on $(D_i \setminus V_i) \times U$. Using these isomorphisms, $Y_1$ and $Y_2$ can be glued to $\Omega\times U$ to form a 3-fold $Z = Y_1 \cup (\Omega\times U) \cup Y_2$.

The $D_i'$ fit together with $\Omega \times \{0\}$ to give an embedded rational curve $L \subset Z$. The germs of $\C^2$-actions also fit together to give the germ of a free $\C^2$-action on $Z$. The orbits are transverse to $L$ everywhere except at the points $0$ and $\infty$.

\begin{figure}[h]
\psfrag{1}{$Y_1$}
\psfrag{2}{$Y_2$}
\psfrag{O}{$\Omega \times U$}
\psfrag{L}{$L$}
\centering
\includegraphics{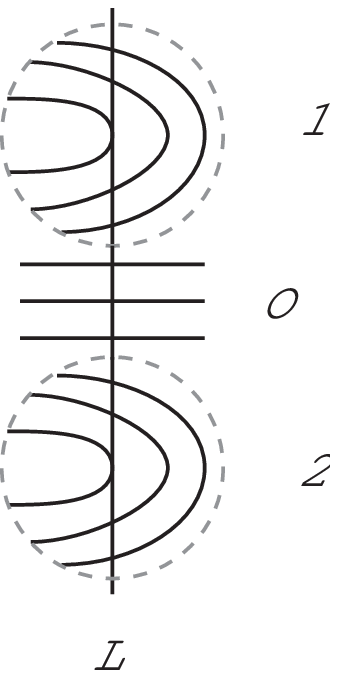}
\emph{\\The twistor space when the $\C^2$-action is free.}

\end{figure}

\begin{lemma}
If $\im \diff\phi_0 \neq \im\diff\psi_\infty$ then $\nu_{L} \cong \mathcal O(1) \oplus \mathcal O(1)$. Otherwise, $\nu_{L} \cong \mathcal O(2) \oplus \mathcal O$.
\end{lemma}

\begin{proof}
Composing the derivative of the $\C^2$-action with the projection map $TZ|_{L} \to \nu_{L}$ gives a map $\C^2 \to \nu_{L}$. The action is free, so the induced map $\Lambda^2 \C^2 \to \Lambda^2 \nu_{L}$ vanishes precisely when the $\C^2$-action fails to be transverse to $L$, i.e., at $0$ and $\infty$. Hence $\Lambda^2 \nu_{L} = \mathcal O(2)$.

Now assume that $\im \diff\phi_0 \neq \im \diff \psi_\infty$. The section of $\nu_{L}$ corresponding to a non-zero point $a\in\im \diff\phi_0$ can only vanish when the $\C^2$-action fails to be transverse. Since $a \notin \im \diff\psi_\infty$, however, it only vanishes at $0$. This section with a single zero gives the short exact sequence
$$
0 \to \mathcal O(1) \to \nu_{L} \to \mathcal O(1) \to 0.
$$
The quotient is $\mathcal O(1)$ since $\Lambda^2 \nu_{L} = \mathcal O(2)$. There are no non-trivial extensions of this kind so $\nu_{L} = \mathcal O(1) \oplus \mathcal O(1)$.

If $\im \diff\phi_0=\im \diff\psi_\infty$, then the section of $\nu_{L}$ corresponding to $a$ vanishes twice giving $\nu_{L} \cong \mathcal O(2) \oplus \mathcal O$.
\end{proof}

Provided then that $\im \diff \phi_0 \neq \im \diff \psi_\infty$, this procedure gives a complex 3-fold containing a rational curve with normal bundle $\mathcal O(1) \oplus \mathcal O(1)$---a twistor space (missing the real structure)---with the germ of a free $\C^2$-action. 

\subsubsection{The real structure}

The antipodal map, $\gamma \colon z \mapsto -\overline z^{-1}$ in a global coordinate, gives $\P^1$ a real structure. Given an involution $\tau$ and a function   $\phi$, taking $\sigma = \gamma \circ \tau \circ \gamma$ and $\psi = \overline{\phi \circ \gamma}$ makes it possible to give the resulting complex 3-fold a real structure, making it into a bona fide twistor space. 

To see this, first note that $\psi$ is defined on $D_2 = \gamma(D_1)$. The pieces $Y_i$ contain embedded copies $D'_i$ of these domains. The map $D_1\times U \to D_2 \times U$ given by $(x, v) \mapsto (\gamma(x), \overline v)$ descends to give a map $\gamma \colon Y_1 \to Y_2$ which extends the antipodal map $D'_1 \to D'_2$. Over $\Omega \times \U$, defining $\gamma(z,v) = (-\overline z^{-1}, \overline v)$ then extends the antipodal map on $L$ to a real structure on the whole of $Z$. Notice that $\diff\psi = \overline{\diff \phi\circ \diff \gamma}$ and so the original assumption on $\diff\phi_0$ ensures that $\im \diff\phi_0 \neq \im \diff \psi_\infty$. So $L$ is a real twistor line in a real twistor space.

The $\C^2$-action and the real structure are related by  $\gamma(v \cdot p) = \overline v \cdot \gamma(p)$. Hence the induced action of $\R^2 \subset \C^2$ commutes with $\gamma$. This means that the $\R^2$-action on the space of twistor lines in $Z$ preserves the real lines. In other words, given a locally defined holomorphic involution $\tau$ of $\P^1$ fixing a single point and a $\C^2$-valued function $\phi$ with $\phi \circ \tau = -\phi$ and $\diff \phi$ non-zero at the fixed point, the above construction produces an anti-self-dual 4-manifold with the germ of a conformally Killing action of $\R^2$

\subsection{The classification}

When taken together, the previous two sections give a classification of germs of toric anti-self-dual 4-manifolds, taken at a point for which the induced $\C^2$-action on twistor space is free. Of course, knowledge of $M$ up to conformal equivalence determines $Z$ up to real-biholomorphism (a biholomorphism which intertwines the real structures) and so these extra equivalences must be taken into account.

Begin by considering the construction of a twistor space without a real structure. Starting with a projective line $\P^1$ and two fixed points $0, \infty \in \P^1$, a complex manifold is built from holomorphic data $(\tau, \phi)$ defined near $0$ and $(\sigma, \psi)$ defined near $\infty$ as above. Two such twistor spaces are bilholomorphic if and only if their defining data are related by an automorphism of $\P^1$ which fixes $0$ and $\infty$.
 
The construction of a real twistor space starts from $\P^1$ with an antipodal map and with fixed antipodal points $0$ and $\infty$; the twistor space is then built from the holomorphic data of a single pair $(\tau, \phi)$. Two pairs give real-biholomorphic twistor spaces if and only if they are related by an automorphism of $\P^1$ which preserves $0$, $\infty$ and the antipodal map. The only such automorphisms are rotations about the axis through $0$ and $\infty$. This determines an action of $S^1$ on the collection of all pairs $(\tau, \phi)$ whose orbits parametrise real-biholomorphism classes of twistor spaces of the relevant kind.

To state the classification: let $\mathcal S$ denote the set of germs of all toric, anti-self-dual 4-manifolds for which the corresponding $\C^2$-action on twistor space is free, taken up to conformal equivalence; let $\mathcal T$ denote the quotient of the set of pairs $(\tau, \phi)$ by the natural $S^1$-action. The twistor correspondence in this context is then:

\begin{theorem}\label{main theorem}
The sets $\mathcal S$  and $\mathcal T$ are in natural one-to-one correspondence.
\end{theorem}

In \cite{joyce}, Joyce completely classifies the local geometry of toric anti-self-dual 4-manifolds under the additional condition that the orthogonal complements to the $\R^2$-orbits form an integrable distribution. As will become clear later, this corresponds precisely to the condition that the involution $\tau$ extends to all of $L$. The exact relationship between Joyce's classification and that given here is discussed in section \ref{joyce construction}. 

\subsection{Toric anti-self-dual K\"ahler metrics}\label{Kahler case}

This section addresses the question of when the conformal class corresponding to a pair $(\tau, \phi)$ admits an invariant K\"ahler representative. As is described in section \ref{review of twistors}, finding a K\"ahler representative of an anti-self-dual conformal class is equivalent to finding a holomorphic section $s$ of $K^{-1/2}$ which is compatible with the real structure and non-zero on real twistor lines. In the case of an invariant K\"ahler metric, $s$ must also be $\C^2$-invariant. This puts restrictions on the possible choice of involution $\tau$.

To see this, begin by considering the situation when the twistor space is without a real structure. Let $Z$ denote the twistor space built from two pairs $(\tau, \phi)$, $(\sigma, \psi)$, and assume that it admits a $\C^2$-invariant section $s$ of $K^{-1/2}$. Let $\theta = s^{-2}$ be the corresponding meromorphic 3-form and let $\alpha = \theta(\tilde X_1 \wedge \tilde X_2)|_L$. Then $\alpha$ is a meromorphic 1-form on the central twistor line $L$ with a simple zeros at $0$ (the fixed point of $\tau$) and $\infty$ (the fixed point of $\sigma$), and double poles at the zeros of $s$ on $L$. 

Since $K^{-1/2}|_L \cong \mathcal O(2)$, $s$ generically has two simple zeros on $L$ (in the real situation described below this follows from the compatibility of $s$ with the real structure). Assume this to be the case and let $z$ denote a global coordinate on $\P^1$ in which $\alpha$ has simple zeros at $0$ and $\infty$ and double poles at $1$ and $c$. Riemann--Roch shows $\alpha$ must be, up to scale,
$$
\alpha_c = \frac{z \,\diff z}{(z - 1)^2(z - c)^2}.
$$

By the Poincar\'e lemma, near $0$, $\alpha = \diff f$ for a unique holomorphic function $f$ with a double zero at $0$. Since $\tau^* \alpha = \alpha$ and $\tau$ fixes $0$, $\tau^* f= f$, hence $\tau$ flips the sheets of the locally defined branched cover given by $f$. This shows that $\tau$ is uniquely determined by $c$; similarly, $\sigma$ is uniquely determined as well. Denote these involutions by $\tau_c$ and $\sigma_c$.

The next step is to work out how many invariant sections of $K^{-1/2}$ there are over a twistor space built from the correct involutions. The 1-form $\alpha_c$ is not the only form preserved by $\tau_c$. If $\lambda , \mu \in \C$ satisfy $\lambda^{-1}\mu = c$ then $\tau_c$ also preserves the 1-form
$$
\beta = \frac{z \diff z}{(z-\lambda)^2(z-\mu)^2}.
$$
This is because the automorphism $g(z) = \lambda^{-1} z$ relates $\alpha_c$ and $\beta$: up to scale, $g^*\alpha_c = \beta$. Hence, near $0$,  $\beta = \diff (f \circ g)$ and, since switching the sheets  of $f$ is the same as switching the sheets of $f \circ g$, $\tau_c^*\beta= \beta$. Similarly, $\sigma_c^* \beta = \beta$. These are the only 1-forms of the relevant type (i.e., with simple zeros at $0$ and $\infty$ and two double poles) preserved by $\tau_c$. This follows from the fact that for $c \neq c'$, $\tau_c \neq \tau_{c'}$. Note that the set of such $\beta$ is parametrised by $\C^*$ by sending $\beta$ to $\mu$.

Hence the twistor space corresponding to the pairs $(\tau, \phi)$ and $(\sigma, \psi)$ can admit an invariant section of $K^{-1/2}$ only when there exists $c \in \C^*$ such that $\tau = \tau_c$ and $\sigma = \sigma_c$. Moreover, on fixing a coordinate on $L$ in which $\tau$ fixes $0$ and $\sigma$ fixes $\infty$, the space of invariant sections, taken up to scale, is parametrised by the subset of $\C^*$ corresponding to the meromorphic 1-forms $\beta$ obtained on $L$. It remains, of course, to show that all such twistor spaces admit invariant sections of $K^{-1/2}$ and that all such meromorphic 1-forms on $L$ are realised. 

Pick $c \in \C^*$ and functions $\phi$, odd with respect to $\tau_c$, and $\psi$, odd with respect to $\sigma_c$. Let $Z$ denote the corresponding twistor space. Pick also a meromorphic 1-form  $\beta$ on $L$, as above, which is preserved by $\tau_c$ and $\sigma_c$. The aim is to find an invariant nowhere vanishing meromorphic 3-form $\theta$, with polar divisor of the form $P = 2Q$ and with $\theta(\tilde X_1 \wedge \tilde X_2)|_L = \beta$. Then $s = \theta^{-1/2}$ is the required holomorphic section of $K^{-1/2}$.

Begin by considering the situation near a point of tangency. Choose a local coordinate in which $\tau_c$ is the map $z \mapsto -z$, defined on the disc $D$. Since $\beta$ is invariant under $\tau_c$, it is the pull back of a 1-form $\beta'$ on $D$ under the map $z \mapsto z^2$. Since $\beta$ has a simple zero at the origin, $\beta'$ is nowhere-vanishing.

Consider $\beta'$ now as a 1-form on $D \times \C^2$ and let $\theta = \beta'\wedge\diff v_1\wedge \diff v_2$, where $(v_1, v_2)$ are coordinates on $\C^2$; it is a $\C^2$-invariant nowhere-vanishing 3-form. Recall that near the point of tangency, $Z$ is the quotient of $D \times \C^2$ where the quotient map $q \colon D\times \C^2 \to D \times \C^2$ is given by $q(z,v) = (z^2, v+\frac{1}{2}\phi(z))$. Hence $q^*\theta = \beta \wedge (\diff v + \frac{1}{2} \diff \phi) = \beta \wedge \diff v$.

The complex 3-fold $Z$ is built by gluing three pieces: $Y_1$, $Y_2$ and $\Omega \times U$. The above discussion gives nowhere-vanishing invariant 3-forms $\theta_i$ on $Y_i$. Define a 3-form over $\Omega \times U$ by $\theta_\Omega = \beta|_\Omega \wedge \diff v_1 \wedge \diff v_2$. The pieces of $Z$ are glued using the quotient maps $q_i \colon D_i \times U \to Y_i$ restricted to $(D_i \setminus V_i) \times U$. Since $q_i^*\theta_i = \theta_\Omega$ the forms $\theta_1, \theta_2, \theta_\Omega$ fit together under this gluing to give a globally defined nowhere-vanishing, $\C^2$-invariant, meromorphic 3-form $\theta$ with double poles along $\{\lambda,\mu\} \times \C^2$. 

This proves the following.

\begin{theorem}
The twistor space corresponding to the pairs $(\tau, \phi)$, $(\sigma, \psi)$ admits a $\C^2$-invariant section of $K^{-1/2}$ with two distinct zeros on each twistor line if and only if there exists $c \in \C\setminus\{0,1\}$ such that $\tau = \tau_c$ and $\sigma = \sigma_c$. In this case, the space of such sections, taken up to scale, is parametrised by $\C^*$.
\end{theorem}

Next consider the case when $Z$ admits a real structure. When the section $s$ of $K^{-1/2}$ is compatible with this real structure, its zeros on $L$ are antipodal. This means that the poles $\lambda$ and $\mu$ of the 1-form $\beta$ must also be antipodal, i.e., $\lambda = -\overline \mu ^{-1}$. This means that the invariant $c= \lambda \mu^{-1}$ is real and negative. Moreover, for such a $c$, the space of $\tau_c$ invariant meromorphic 1-forms compatible with the real structure is parametrised by $S^1$. 

In order to check that all such $c$ and $\beta$ actually arise it is necessary to check that the section $s$ constructed above is compatible with the real structure. This is more-or-less immediate from the definition of the anti-holomorphic involution and the fact that $\beta$ is compatible with the antipodal map. This gives the following classification of toric, anti-self-dual, K\"ahler surfaces:

\begin{theorem}\label{conformally Kahler}
The conformal class $(\tau, \phi)$ admits a K\"ahler representative if and only if $\tau = \tau_c$ for some real number $c <0$. In this case the space of distinct K\"ahler representatives, taken up to scale, is parametrised by $S^1$.
\end{theorem}

An interesting special case of this situation is $c = -1$. The compatible forms $\beta$ have no residues, so the function $f$ satisfying $\beta = \diff f$ is a globally defined meromorphic function $f \colon L \to \P^1$. Since $\tau_{-1}$ flips the sheets of this branched cover, it extends to the whole of $L$. In a global coordinate $z$ on $L$ in which the anitpodal map is $\gamma(z) = -\overline z^{-1}$, the involution is $\tau_{-1}(z) = -z$.

As will be seen after the metric has been computed explicitly (in section \ref{metric formulae}), this means that the distribution of orthogonal complements to the orbits in $M$ is integrable. This is precisely the case considered by Joyce in \cite{joyce}, where the $S^1$ of K\"ahler representatives are also constructed.

Finally, the form of $\tau$ determines whether or not the orbits in $M$ are Lagrangian.

\begin{lemma}\label{lagrangian orbits}
A K\"ahler representative for the conformal class $(\tau_c, \phi)$ has Lagrangian orbits if and only if $c = -1$.
\end{lemma}

\begin{proof}
Let $\omega$ denote the K\"ahler metric. The function $\omega(X_1, X_2)$ is a Hamiltonian for $[X_1, X_2]=0$ and hence is constant. From section \ref{review of twistors}, at $p \in M$, $\omega(X_1, X_2) = 2\pi i \res \beta$ (up to sign). Hence the orbits are Lagrangian if and only if $\beta$ has no residues. This only happens for $c=-1$.
\end{proof}

\subsection{An example: the Fubini--Study metric}

The projective plane $\P^2$ with the Fubini--Study metric is toric and \emph{self-dual}. This difference, which is just one of orientation, does not affect the previous discussion of holomorphic data. This section computes the corresponding involution $\tau$ and odd function $\phi$. 

The twistor space of $\P^2$ is the flag manifold 
$$\
\F = \{(p,l) : l \subset \P^2\ \text{a line,}\ p\in l\}.
$$
The twistor lines in $\F$ correspond to pairs $(q, k)$ where $q \in \P^2$ and $k\subset \P^2$ is a line with, this time, $q\notin k$; the line $L \subset \F$ corresponding to $(q,k)$ is the set
$$
L = \{(p,l) : p \in k,\ l = \overline{pq}\}.
$$
The real structure is given by $(p,l) \mapsto (l^{\perp}, p^{\perp})$.
Real twistor lines correspond to pairs $(k^{\perp}, k)$; in this case the line is naturally identified with $k$ itself. The twistor projection to $\P^2$ is given by $(p, l) \mapsto p^{\perp}\cap l$ where $p^{\perp}\cap l$ is the orthogonal complement of $p$ in $l$.

The Killing fields $X_i$ are induced by a $T^2$-action on $\P^2$; the holomorphic fields $\tilde X_i$ are thus induced by a $(\C^*)^2$-action on $\F$. Hence, given a twistor line $L$, there is a (a priori only locally defined) function $\lambda \colon L \to (\C^*)^2$ such that $\lambda(z)\cdot z \in L$. The function $\phi$ is determined by $\lambda = e^{i\phi}$ and $\phi(0) = 0 = \phi(\infty)$. The odd condition $\phi \circ \tau = -\phi$ translates into $\lambda \circ \tau = \lambda^{-1}$, whilst the real-compatibility condition $\phi\circ \gamma = \overline\phi$ becomes $\lambda \circ \gamma = \overline{\lambda}^{-1}$.

The $(\C^*)^2$-action induced on $\F$ by considering it as the twistor space of the toric, self-dual, manifold $\P^2$ is the same as the obvious $(\C^*)^2$-action. Using this it is straightforward to compute $\lambda$.

Let $L$ denote a real twistor line $(k^{\perp},k)$ in $\F$. Provided $k$ is not one of the three toric divisors, the action on $\F$ is free near $L$. Choose coordinates $\zeta_1, \zeta_2$ on $(\C^*)^2$ and denote the corresponding components of $\lambda$ by $\lambda_1, \lambda_2$. Since the $(\C^*)^2$-orbits exist globally on $\F$, the $\lambda_i$ are global meromorphic functions. Moreover, generic orbits meet $L$ twice, so the $\lambda_i$ have degree 2. Notice that this implies that $\tau$ extends to all of $L$. It remains to find the poles and zeros of the $\lambda_i$.

Let $l_1, l_2, l_3$ be the toric divisors in $\P^2$ with stabilizers $\zeta_1=1, \zeta_2=1, \zeta_1= \zeta_2$ respectively. There are six distinguished points on $L$, namely those corresponding to the points $p_i = k \cap l_i$ and $\gamma(p_i)$. Since the $p_i$ have stabilizers under the action, it follows that $\lambda_1$ has a pole at $p_1$, $\lambda_2$ has a pole at $p_2$ and both have a zero at $p_3$. By the reality condition, it follows that $\lambda_1$ has a zero at $\gamma(p_1)$, $\lambda_2$ has a zero at $\gamma(p_2)$ and both have a pole at $\gamma(p_3)$. Since they each have degree 2 all poles and zeros are simple. 

It is now possible to write the $\lambda_i$ down explicitly. Choose a global coordinate $z$ on $k$ in which $\tau(z) =-z$, $\gamma(z) = -\overline z^{-1}$; then, for some $c_i$,
$$
\lambda_i(z)
= c_i \frac{(z + \overline p_i^{-1})(z - p_3)}{(z - p_i)(z + \overline p_3^{-1})}.
$$
The conditions $\lambda_i(0) = 1 = \lambda_i(\infty),$ $\lambda \circ \tau = - \lambda$ and $\lambda \circ \gamma = \overline \lambda ^{-1}$ force $c_i=1$ and $|p_i| = 1$ for $i=1,2,3$. Rotating so that $p_3 = 1$ uniquely determines the coordinate $z$ and gives: 
$$
\lambda_i(z)
=
\frac{(z + p_i)(z-1)}{(z-p_i)(z+1)}
\quad
\text{for}\ p_1, p_2 \in S^1\setminus{\pm 1}.
$$
The parameters $p_1, p_2$ correspond to the fact that there is a two parameter family of real lines modulo the $T^2$-action.

\section{Metric formulae}
\label{metric formulae}

This section explains how, in principle at least, the metric can be computed explicitly from the corresponding holomorphic data. For notational convenience, from now on $\tau$ will denote the involution defined near $0$ and $\infty$ in $L$ and $\phi$ the odd function (even in the purely complex case where $\tau$ and $\phi$ are not invariant under the antipodal map).

\subsection{Describing the space of lines}

In order to recover the anti-self-dual conformal class from its twistor space $Z$, the first step is to describe the space of twistor lines in $Z$.

To begin with, concentrate on the purely complex situation; the real involution will be re-introduced later. Let $(\tau, \phi)$ be holomorphic data defined near $0$ and $\infty$ determining a twistor space $Z$. Let $L'$ be a line in $Z$ near to the central line $L$. There is a $\C^2$-orbit which is tangent to $L'$ at a point near to $0\in L$; this orbit meets $L$ in a pair of points $a$ and $\tau(a)$ near $0$. Likewise, there is another pair of distinguished points $b, \tau(b) \in L$ near $\infty$.

\begin{figure}[h]
\psfrag{a}{$a$}
\psfrag{ta}{$\tau(a)$}
\psfrag{b}{$b$}
\psfrag{tb}{$\tau(b)$}
\psfrag{L'}{$L'$}
\centering
\includegraphics{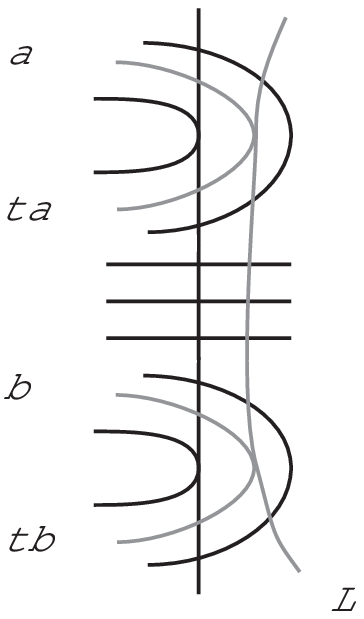}
\emph{\\A line near to the central line.}

\end{figure}

The $\C^2$-action on $Z$ defines a $\C^2$-action on the space of lines. Other lines near $L$ and in the same orbit as $L'$ will lead to the same four points of $L$. Thus the quotient of the space of lines by the $\C^2$-action is identified with pairs $(a, b)$, or more precisely, the quotient of these by the $\Z_2\times \Z_2$-action generated by $\tau$ acting on each of the two components of its domain. To fix notation, let $r$ denote a coordinate on the quotient of the domain of $\tau$ near $0$ by $\tau$ itself; let $s$ be a similar coordinate for the quotient near $\infty$. The quotient of the space of lines by the $\C^2$-action is identified with the pairs $(r, s)$.

For a fixed pair $(r,s)$, it is possible to describe the corresponding lines as follows. Over the middle portion $\Omega \times U$ of $Z$ the $\C^2$-action takes each point of $L'$ to a unique point of $L$; hence $L'$ is the graph of a holomorphic map $f \colon \Omega \to \C^2$. Near $0$, however, the $\C^2$-orbits meet $L$ twice. To extend this description near $0 \in L$ it is necessary to make a cut from $a$ to $\tau(a)$. Then there is a well defined extension of $f$ which jumps 
\begin{equation}\label{jump}
f(z) \to f(\tau(z)) + \phi(z)
\end{equation}
as $z$ moves over the the cut. Similarly, after cutting $L$ from $b$ to $\tau(b)$, $f$ extends near $\infty$ and satisfies the analogous jump condition over the cut there.

This ``jumping condition'' can be interpreted in terms of prescribing the \v{C}ech coboundary of a cochain on a certain elliptic curve. The most natural approach is to construct the elliptic curve by taking two copies of the cut plane and gluing them using $\tau$. This method is described below; it leads, however, to some techincally awkward considerations. In the special case that $\tau$ extends to the whole twistor line, it is possible to work instead with the elliptic curve which branches over $L$ with branch points at $a, \tau(a), b \tau(b)$---that is the elliptic curve obtained by gluing two copies of the cut plane directly and \emph{not} via $\tau$. This is technically more straightforward, but contains the essential ideas present in the general case. For this reason it is described first.

\subsection{The space of lines when $\tau$ extends}\label{when tau extends}

Assume that $\tau$ extends to the whole of $L$ and choose a coordinate on $L$ in which $\tau(z) = -z$. Let $E \to L$ denote the elliptic curve branched over $L$ at the points $\pm a,\pm b$.  The elliptic curve has two natural involutions: the first, $\rho$, comes from its description as a branched cover; the second is the lift of $\tau$, which is also denoted $\tau$. In the model of $E$ as pairs $(z,w)$ with $w^2 = (z^2 - a^2)(z^2-b^2)$ these involutions are given by $\rho(z,w) =(z, -w)$, $\tau(z,w) = (-z,w)$. Note that $\tau$ and $\rho$ commute. 

Cutting $L$ from $a$ to $-a$ and from $b$ to $-b$ determines two open annuli $U, V \subset E$ which form a cover; moreover they can be chosen so that $V = \rho(U)$. Pulling back $\phi$ to $E$ defines a \v{C}ech 1-cochain $\hat \phi$ with respect to this cover. It satisfies $\hat \phi \circ \rho = \hat \phi$ and $\hat \phi \circ \tau = - \hat \phi$.

Using this set up it is possible to describe functions which satisfy the jumping condition (\ref{jump}) as \v{C}ech 0-cochains with coboundary $\hat \phi$.

\begin{lemma}
There exists a 0-cochain $(f_U, f_V)$ with coboundary $\hat \phi$. It is unique up to an additive constant, and satisfies $f_V=f_U \circ\tau\circ\rho$
\end{lemma}

\begin{proof}
The open sets $U$, $V$ and $U\cap V$ are all (disjoint unions of) annuli and so Stein; hence $\{ U, V\}$ gives an acyclic cover of $E$ and so can be used to compute the cohomology groups $H^*(\mathcal O)$. Since $H^0 = \C = H^1$ this implies that $(f_U, f_V)$ exists if and only if $\hat \phi$ satisfies a single linear condition and, if this condition is satisfied, $(f_U, f_V)$ is unique up to the addition of a constant. 

To formulate the condition, let  $\omega$ be a regular differential on $E$ and $C$ a contour in $E$ with two components, one going round each cut. If $(f_U, f_V)$ exists,
$$
\int_C \hat \phi \, \omega
=
\int_C f_U\omega - \int_C f_V \omega
=0
$$
since $f_U \omega$ is holomorphic over $U$ whilst $f_V\omega$ is holomorphic over $V$. So $\int \hat\phi\,\omega=0$ is a necessary and hence sufficient condition for the existence of $(f_U, f_V)$.

For the 1-cochain $\hat \phi$ considered here, this condition is satisfied since $\hat \phi$ is $\tau$-odd, whilst $\omega$ is $\tau$-even as can be seen, for example, from the formula
$$
\omega = \frac{\diff z}{(z^2 - a^2)^{1/2}(z^2 - b^2)^{1/2}}.
$$
Hence $(f_U, f_V)$ exists and is unique up to a constant.

To verify that $f_V = f_U\circ \tau \circ\rho$, define $F_U = f_V\circ \tau\circ \rho $ and $F_V = f_U \circ \tau \circ \rho$. Then, on $U \cap V$, $F_U - F_V = - \hat \phi\circ\tau \circ\rho = \hat \phi$. By uniqueness of $(f_U, f_V)$, there is a constant $c$ such that $f_U = f_V\circ\tau\circ\rho + c$ and $f_V=f_U \circ\tau\circ\rho + c$. Hence $f_U = f_U + 2c$ and so $c =0$.
\end{proof}

This shows that the solutions to the ``jumping problem'' (\ref{jump}) are exactly the 0-cochains $(f_U, f_V)$. It also shows that they always exist and are parametrised by the constants $\C^2$ (remembering that $\phi$ is $\C^2$-valued). Of course, this existence statement is guaranteed by the general theory (there are lines in $Z$ nearby to $L$ by Kodaira's deformation theory); in this particular situation, however, it also follows from the elementary \v{C}ech cohomological argument above.

More importantly, it provides a concrete desription of the space of lines: explicitly finding 0-cochains with presecribed coboundary is a classical problem which is solved by contour integration yielding, in this case at least, classical integral formulae. This is demonstrated by the following result, in which $\alpha(p)$ is a meromorphic 1-form on $E$ with poles at $p$ and $\rho(p)$ with residues $1$ and $-1$ respectively and $C$ is a contour with two components, one round each cut. 

\begin{lemma}
$f_U(p) = c + \frac{1}{4\pi i} \int_C \phi\, \alpha(p)$, for some constant $c$ where $\gamma$ is chosen so that $p$ and $\rho(p)$ lie on opposite sides
\end{lemma}

\begin{proof}
Let $p\in U$. By Cauchy's residue theorem applied to $f_U \alpha$ over $U$ and $f_V \alpha$ over $V$, $\frac{1}{2\pi i}\int_C \phi\,\alpha(p) = f_U(p) -f_V(\rho(p))$.

On the other hand, there is exists a constant $c$ such that $f_U+ f_V\circ \rho = c$. To see this, define $F_U = -f_V\circ \rho$ and $F_V = -f_U\circ\rho$. Then, over $U \cap V$, $F_U - F_V = \hat \phi \circ \rho = \hat \phi$. By the uniquness of $(f_U, f_V)$ this gives that $f_U(p) + f_V ( \rho(p)) =c$

Adding these two equations gives the result.
\end{proof}

Using this formula, it is possible to explicitly compute the conformal class in the coordinates $(r, s, v)$. Recall from section \ref{review of twistors} that it suffices to find the null cone and that the tangent vector $(r',s',v')$ at $(r,s,v)$ is null if and only if the corresponding section of the normal bundle of the line $(r,s,v)$ has a zero. The line determined by $(r,s,v)$ is the ``graph'' of a function $f$ satisfying the jump condition (\ref{jump}) over cuts between $\pm\sqrt r$ and between $\pm \sqrt s$. The section of the normal bundle corresponding to $(r',s',v')$ is given by $r' \del f/ \del r + s' \del f/\del s + v'$. This can be computed explicitly by differentiating the formula from the previous lemma.

\begin{lemma}
$$
\frac{\del f}{\del r}(p)
=
A\left(\frac{p^2 - s}{p^2 - r}\right)^{1/2},
$$
where $A \in \C^2$ is given by 
$$
A 
=
-\frac{1}{8\pi i}
\int_C
\frac{ z \phi(z)\,\diff z}
{(z^2 - r)^{3/2}(z^2- s)^{1/2}}.
$$
\end{lemma}
\begin{proof}
Differentiating the integral formula for $f$ gives
$$
\frac{\del f}{\del r}(p) 
= 
\frac{1}{4\pi i }\int_\gamma \phi \frac{\del \alpha}{\del r}(p).
$$
The meromorphic $1$-form $\alpha(p)$ is given by
$$
\alpha(p) = 
\frac{\diff z}{z - p}
\sqrt{
\frac
{(p^2 - r)(p^2 - s)}
{(z^2 - r)(z^2 - s)}
}.
$$
Differentiation with respect to $r$ gives
$$
\frac{\del \alpha}{\del r}(p) = 
\theta(z) (p+z)
\left( \frac{p^2 - s}{p^2 - r} \right)^{1/2}
$$
where 
$$
\theta =  
\frac{ - \diff z}
{2(z^2 - r)^{3/2}(z^2 - s)^{1/2}},
$$
The $p \theta$ term does not contribute to $\int \phi\frac{\del\alpha}{\del r}$ since $\phi\theta$ is odd. Hence,
$$
\frac{\del f}{\del r}(p) 
=
A\left (\frac{p^2 - s}{p^2 - r} \right)^{1/2},
$$
as claimed.
\end{proof}
There is, of course, a similar formula for $\del f/ \del s$ with the roles of $r$ and $s$ swapped; that is
$$
\frac{\del f}{\del s}(p)
=
B \left(\frac{p^2 - r}{p^2-s}\right)^{1/2},
$$
where 
$$
B = 
-\frac{1}{8\pi i}
\int_C
\frac{ z \phi(z)\,\diff z}
{(z^2 - r)^{1/2}(z^2- s)^{3/2}}.
$$

The null cone is then the set of $(r',s',v')$ for which the following equations have a simultaneous solution for some $p$:
$$
r' A_i \left(\frac{p^2- s}{p^2 - r}\right)^{1/2}
+
s' B_i \left(\frac{p^2- r}{p^2 - s}\right)^{1/2}
+
v'_i
=0,
$$
where $A_i, B_i, v'_i$ are the components of $A, B, v$. Eliminating $p$ gives the following formula for the conformal class.

\begin{theorem}\label{metric when tau extends}
Let $(\tau, \phi)$ determine a twistor space in which $\tau$ extends to the whole of the central twistor line. Then the conformal class on the space of lines is given, in the coordinates $(r,s,v)$ described above, by
\begin{equation}\label{so metric}
\diff r\diff s + 
\frac{
(A_2\diff v_1 - A_1\diff v_2)(B_2\diff v_1 - B_1\diff v_2)}
{(A_2 B_1 - A_1B_2)^2}
\end{equation}
where $A_i, B_i$ are the components of the contour integrals $A, B$ from above.
\end{theorem}

It should be noted that we have been slightly cavalier in our approach to the square roots involved in the definition of $A$ and $B$. In fact, $A$ and $B$ are defined only up to sign, with $A$ behaving like $s^{-1/2}$ and $B$ like $s^{-3/2}$ as $s$ approaches infinity. The sign ambiguity disappears in the products $AB$ which appear in the metric formula.

\subsection{Surface orthogonal actions and Joyce's construction}\label{joyce construction}

The metric formula (\ref{so metric}) shows that when $\tau$ extends to the whole of the central twistor line, the orthogonal distribution to the $\C^2$-action is integrable. Such an action is said to be \emph{surface-orthogonal}. In \cite{joyce}, Joyce classifies the local geometry of such anti-self-dual conformal classes. In \cite{calderbank-pedersen} Calderbank and Pedersen give the following description of Joyce's work.

Generically, surface-orthogonal anti-self-dual conformal classes are determined by a pair of axially-symmetric harmonic functions defined on an open set in ${\R}^{3}$. So we take as input a vector-valued function $F=F(x,y)$ on an open set in the upper-half plane $y>0$ satisfying the equation
\begin{equation}\label{ash}
F_{xx}+F_{yy}+ \frac{1}{y} F_{y}= 0.
\end{equation}
Here $x,y$ are real co-ordinates and $F$ takes values in $\R^{2}$. Set $P= - y F_{x}$ and $Q= y F_{y}$ so that the derivatives of $P,Q$ satisfy the system of linear equations, the \emph{Joyce equations}:
$$
P_{x}= Q_{y}\ \ ,\ \  P_{y}+Q_{x}= y^{-1} P. 
$$
Conversely, any solution of these equations arises (at least locally) from some $F$. Now write $P_{i}, Q_{i}$ for the components of $P,Q$. Then Joyce's metric has the form
$$
\frac{\diff x^{2}+\diff y^{2}}{y^{2}} 
+ 
\frac{
(P_{2} \diff u_{1}-P_{1} \diff u_{2})^{2}
+
(Q_{2} \diff u_{1}- Q_{1} \diff u_{2})^{2}}
{(P_{1} Q_{2} - Q_{1} P_{2})^{2}}. $$
This is an anti-self-dual Riemannian metric on a real 4-dimensional manifold; the anti-self-duality condition is a consequence of the Joyce equations.

We now complexify Joyce's construction to produce an anti-self-dual conformal structure on a 4-dimensional complex manifold. Clearly we can just allow $x,y$ to be complex variables and $F$ to be a $\C^{2}$-valued solution of (\ref{ash}). Change variables by writing $\zeta=x+iy, \xi=x-iy$. Then the equation satisfied by $F$
becomes
\begin{equation}\label{zeta-xi equation}
F_{\zeta \xi}= \frac{1}{2(\zeta-\xi)}(F_{\zeta}-F_{\xi}).
\end{equation}
Put $A=F_{\zeta}, B=F_{\xi}$. The Joyce metric is conformal to
$$
\diff \zeta \diff \xi +  H_{A,B}
$$
where 
$$
H_{A,B}
= 
\frac
{(A_{2} \diff u_{1}- A_{1} \diff u_{2})
(B_{2} \diff u_{1} - B_{1} \diff u_{2})}
{(A_{2}
B_{1}- A_{1} B_{2})^{2}}.
$$

This is precisely the same expression obtained in Theorem \ref{metric when tau extends} with a metric $\diff r \diff s +H_{A,B}$ where $A_{i}, B_{i}$ were given by contour integrals. To match things up further, observe that for any fixed $z$ the function of $r,s$ given by
$$
g(r,s)= \frac{1}{(z^{2}-r)^{1/2}(z^{2}-s)^{1/2}},
$$
satisfies the equation
\begin{equation}\label{r-s equation}
g_{rs}=\frac{1}{2(r-s)} (g_{r}- g_{s}).
\end{equation}
If we write
$$
G(r,s) = 
- \frac{1}{4\pi i} \int_{C} 
\frac{z\phi(z)}{(z^{2}-r)^{1/2}(z^{2}-s)^{1/2}}\diff z,
$$
where $\phi$ and  the contour $C$ are as considered in the previous section then $A=G_{r}, B=G_{s}$ and, differentiating under the integral sign, $G$ satisfies the same equation $G_{rs}=\frac{1}{2(r-s)} (G_{r}-G_{s})$ as $F$, when we replace the co-ordinates $r,s$ by $\zeta, \xi$. Notice that, similar to the discussion at the end of the previous section, $G$ is only defined up to a sign and behaves like $s^{-1/2}$ as $s$ tends to $\infty$.

On closer inspection, however, the correspondence between the two points of view is less straightforward. The key is to understand the transformation of the equation (\ref{zeta-xi equation}) under Mobius maps. In general, let $f(z)= \frac{az+b}{cz+d}$ be a Mobius map and consider two pairs of variables $r=f(\zeta), s=f(\xi)$. Suppose a function $G(r,s)$ satisfies equation (\ref{r-s equation}) and, changing variables, put $\alpha=G_{\zeta}, \beta=G_{\xi}$. Then we have 
$$ 
\alpha_{\xi}
= 
r_{\zeta} s_{\xi}  G_{rs}
=
\frac{1}{2(r-s)} \left( s_{\xi} \alpha - r_{\zeta} \beta\right), 
$$
where we write $r_{\zeta}$ for the derivative $\diff r/\diff \zeta = f'(\zeta)$. Clearly we also have the identity $\alpha_{\xi}=\beta_{\zeta}$. We now seek a function $\lambda(\zeta,\xi)$
such that $\alpha'=\lambda \alpha$ and $\beta'=\lambda^{-1} \beta$ satisfy
\begin{equation}\label{alpha-beta equation}
\alpha'_{\xi}
= 
\frac{1}{2(\zeta-\xi)} (\alpha'-\beta')
= 
\beta'_{\zeta}.
\end{equation}
If we can do this then there is a function $F(\zeta,\xi)$ such that $\alpha'=F_{\zeta}, \beta'=F_{\xi}$ and $F$ satisfies equation (\ref{zeta-xi equation}). 

To find $\lambda$ we expand the first equation in (\ref{alpha-beta equation}) to get
$$ 
\lambda \alpha
\left( 
2 \lambda^{-1} \lambda_{\xi} 
+ 
\frac{s_{\xi}}{r-s} 
- 
\frac{1}{\zeta-\xi}
\right)
-
\beta\left(
\frac{\lambda r_{\zeta}}{r-s}- \frac{\lambda^{-1}}{\zeta-\xi}\right)=0. 
$$
We seek a solution in which the co-efficients of $\alpha,\beta$ both vanish. The $\beta$ coefficient gives
$$ 
\lambda^{2}= \frac{(r-s)}{ r_{\zeta} (\zeta-\xi)}. 
$$
With this definition of $\lambda$ we have
$$ 
2 \lambda^{-1} \lambda_{\xi}= - \frac{s_{\xi}}{r-s} - \frac{1}{\zeta-\xi}
$$
so the $\alpha$ coefficient also vanishes and the first equation in (\ref{alpha-beta equation}) is satisfied. Symmetrically, we can solve the second equation if
$$ 
\lambda^{-2}= \frac{(r-s)}{s_{\xi}(\zeta-\xi)}. 
$$

So far we have not used the fact the the variables are related by a Mobius map. This enters now in the compatability between the two formulae for $\lambda$. For any Mobius map $f$ we have, as one easily checks, the identity
$$  
f'(\zeta) f'(\xi)= \left(\frac{ f(\zeta)-f(\xi)}{\zeta-\xi}\right)^{2}. 
$$
Thus the two equations above are compatible and we can find $\alpha',\beta'$ and hence, implicitly, the function $F$. 

In terms of our metrics on 4-dimensional complex manifolds, this means that if we start with a metric 
$$
\diff r \diff s+H_{G_{r}, G_{s}}
$$
determined by a function $G(r,s)$ which satisfies (\ref{r-s equation}) then we can change coordinates on the base of the fibration to $\zeta, \xi$ and the metric is conformal to 
$$
\diff \zeta \diff \xi+H_{F_{\zeta},F_{\xi}}
$$ 
determined by a function $F(\zeta,\xi)$ where $F$ is a solution of (\ref{zeta-xi equation}). This is because 
$$ 
\diff r \diff s+H_{G_{r}, G_{s}}
= 
r_{\zeta}s_{\xi}\left(
\diff \zeta \diff \xi 
+ 
H_{ r_{\zeta} G_{r}, s_{\xi} G_{s}}
\right)
$$
and 
$$
H_{r_{\zeta} G_{r}, s_{\xi} G_{s}}
= 
H_{G_{\zeta}, G_{\xi}}
=
H_{\lambda G_{\zeta}, \lambda^{-1} G_{\xi}}
= 
H_{F_{\zeta}, F_{\xi}}. 
$$

We now implement this procedure in our case, so $r$ lies in a neighborhood of $0$ and $s$ in a neighbourhood of $\infty$. We take 
$$  
r=\frac{\zeta-i}{\zeta+i},\ \  s=\frac{\xi-i}{\xi+i}, 
$$
so $\zeta, \xi$ lie in neighborhoods of $\pm i$ respectively. Then
$$ 
\frac{(r-s)}{s_{\xi}(\zeta-\xi)} = \frac{\xi+i}{\zeta+i}, 
$$
so $\lambda= \left(\frac{\zeta+i}{\xi+i}\right)^{1/2}$. The function $F(\zeta,\xi)$ is determined (up to an arbitrary constant) by the equations
\begin{equation}\label{F' equations}  
F_{\zeta}
= 
\left(\frac{\zeta+i}{\xi+i}\right)^{1/2} G_{\zeta},\quad
F_{\xi}
= 
\left(\frac{\xi+i}{\zeta+i}\right)^{1/2}G_{\xi}. 
\end{equation}
Here we should consider the indeterminacy involved in the square roots. We can assume that $\zeta$ does not take on the value $-i$, so the only branching is at $\xi=-i$, which corresponds to $s=\infty$. Now recall that the function
$G(r,s)$ is itself not defined across $s=\infty$, but behaves like $s^{-1/2}$. Then one easily checks that the products involved in (\ref{F' equations}) are well-defined, so $F$ is a genuine holomorphic function of $\zeta,\xi$ around $\zeta=i, \xi=-i$. 

Finally, we can go back and consider the case of data compatible with real structures. In the $r,s$ variables the real points are given by $s=\overline{r}^{-1}$, which goes over to $\xi=\overline{\zeta}$. Recall that the reality condition for our original function $\phi$ is $\phi(-\overline{z}^{-1})= \overline{\phi(z)}$. Changing the variable $z$ in the contour integral to $\overline{z}^{-1}$ shows that the derivatives of $G$ satisfy the condition
$$
G_{\overline{r}}= \frac{1}{\sqrt{r\overline{r}}} \overline{G_{r}}. 
$$
Then a few lines of calculation shows that the derivatives of $F$ obey
$$ 
F_{\overline{\zeta}}= \overline{F_{\zeta}},
$$
so the imaginary part of $F$ is constant and without loss of generality $F$ is real. 

To sum up, we have seen that our twistor analysis yields an alternative proof of Joyce's description  of surface-orthogonal solutions (in the generic case). What is missing is an explicit formula for Joyce's function $F$ in terms of our holomorphic data. It is a classical fact, explained to us by Hitchin,
that any local solution $F(\zeta, \xi)$ of (\ref{zeta-xi equation}) can be expressed as a contour integral
\begin{equation}\label{integral rep for ash}
F(\zeta,\xi)
= 
\int_{\Gamma} \frac{\Psi(u)}{(u-\zeta)^{1/2}(u-\xi)^{1/2}} \diff u, 
\end{equation}
where the the square root is defined by cutting the $u$-plane from $\zeta$ to $\xi$ and $\Gamma$ is a contour encircling $\zeta, \xi$. If we write $u=v^{2}$ we have a representation very much like (\ref{integral rep for ash}):
$$
F(\zeta,\xi)
= 
\int_\Gamma\frac{v\Psi(v^{2})}{(v^{2}-\zeta)^{1/2}(v^{2}-\xi)^{1/2}} dv, 
$$
but with a different contour. It is tempting to try to relate
the holomorphic data $\Psi$ and $\phi$ but we have not yet succeeded in doing this.

\subsection{The space of lines for general $\tau$}

Now return to the case of general $\tau$. In this situation, the jumping condition (\ref{jump}) is most naturally interpreted on the elliptic curve $E$ produced by gluing two copies of the cut plane via the involution $\tau$. This curve is, in a natural way, a branched double cover of the Riemann sphere, but---unlike in the previous case---this is not the map which extends the inclusion on each piece of $E$.

\begin{proposition}
There is a natural branched double cover $E \to \P^1$.
\end{proposition}

\begin{proof}
Let $D$ be the disc about $0$ on which $\tau$ is defined, and let $A$ be the double branched cover of $D$, branched at $a$, $\tau(a)$. Then $A$ is an annulus and each component of the boundary of $A$ is mapped bijectively to the boundary of $D$. By definition, $A$ carries a holomorphic involution $\rho$, i.e., $D = A/ \rho$. 

The involution $\tau$ lifts to an involution $\tilde \tau$ of $A$ which maps each boundary component of $A$ to itself and commutes with $\rho$. This can be seen explicitly. Take a local coordinate $w$ on $D$ in which $\tau(w) = -w$ and suppose that $a$ is the point $w=2$. Then $A = \{\tilde w \in \C : 1/2 < |\tilde w|<2\}$ and $\rho(\tilde w) =  \tilde w ^{-1}$ whilst the covering map is given by $w = \tilde w + \tilde w^{-1}$. The involution $\tilde \tau$ is given by $\tilde \tau (\tilde w) = -\tilde w$. 

Now $\rho' = \rho \circ\tilde \tau$ is another involution of $A$ which exchanges the boundary components. The quotient $A/ \rho'$ is again biholomorphic to a disc $D'$, the boundary of which is canonically identified with the boundary of $D$. Cutting out $D$ from $L$ and gluing in $D'$ using this identification, and doing the same near $\infty$, gives a new Riemann surface $\P^1_{r,s}$.

Of course, abstractly, $\P^1_{r,s}$ is just another copy of the Riemann sphere; the point is that now the inclusion maps to $L$ on the two halves of the elliptic curve $E$ extend to give a branched double cover $E \to \P^1_{r,s}$. 
\end{proof}

In fact, $\P^1_{r,s}$ can be canonically identified with a line in the twistor space corresponding to the parameters $r,s$. That is, the obvious identification of the big open set $\Omega \subset L$ with an open set in the other twistor line, given by the $\C^2$-action, extends if $\Omega$ is regarded as a subset of $\P^1_{r,s}$.

The earlier case, in which $\tau$ is globally defined, is seen to be special in the following way. The elliptic curve made in that case has two holomorphic involutions, in fact an action of $\Z_2 \times \Z_2$, corresponding in the model $w^2 = (z^2-r)(z^2 - s)$ to $z\mapsto \pm z,\ w\mapsto\pm w$. Thus there are different descriptions of $E$ as a branched cover. In the general case there is only one such description, and this is \emph{not} the one used previously.

With this in hand, it is now possible to show how the twistor lines  can be described via solutions to the problem of prescribing the \v{C}ech coboundary of a cochain on $E$. 

Denote by $\rho'$ the involution of $E$ such that $E/ \rho' = \P^1_{r,s}$. Let $U$ and $V$ denote the open cover of $E$ coming from its construction as two copies of $L$ glued along cuts via $\tau$. Note that $\rho'$ exchanges $U$ and $V$. Restrict $\phi$ to the disc $D$ and then lift to the annulus $A$. This gives a function $\tilde \phi$ which, regarding $A$ as a subset of $E$, is a function on one component of $U\cap V$. Do the same to define $\tilde \phi$ on the other component of $U \cap V$; this gives a \v{C}ech 1-cochain which is odd with respect to $\rho'$.

Suppose now that $f$ is a $\C^2$-valued function on $L$ which satisfies the ``jumping condition'' over the cuts. Write $f_U$ for $f$ thought of as a function on $U\subset E$ and $f_V$ for $f$ thought of as a function on $V\subset E$. Then $\tilde \phi = f_U - f_V$ and the 1-cochain is seen to be a coboundary. Conversely, finding a 0-cochain with $f_V=f_U\circ \rho'$ and coboundary $\tilde \phi$ gives a solution to the original problem.

As before, the prescribed coboundary problem has a unique solution, modulo constants. It remains to be checked that the functions satisfy $f_V = f_U \circ \rho'$ but this follows from the fact that $\tilde \phi$ is odd with respect to $\rho'$.

Moreover, the functions $f_U$ can again be found via contour integration. In the current situation this is done as follows. Let $q_*$ be a fixed base point in $\P^1_{r,s}$ and $q$ any other point. Write $\alpha_q$ for the unique meromorphic 1-form on $\P^1_{r,s}$ with simple poles at $q$ and $q_*$ of residues $1$ and $-1$ respectively. Write $\tilde \alpha_q$ for the pull back of $\alpha_q$ to $E$ via the branched cover and denote by $\tilde q \in U$ the inverse image of $q$.
 
\begin{lemma}\label{formula for f}
Up to an additive constant, 
$$
f(\tilde q) = \frac{1}{4\pi i}\int_{\gamma}  
\tilde\phi\,\tilde\alpha_q.
$$
\end{lemma}
\begin{proof}
Write $\tilde\phi = f-f\circ\rho'$. Use Cauchy's residue theorem over $U$ to evaluate the integral of $f\tilde\alpha_q$ and over $V$ to evaluate the integral of $(f \circ \rho') \tilde \alpha_q$.
\end{proof}

The next step is to compute the conformal class in the coordinates $(r,s,v)$. This amounts to determining the null cone, which completely determines the conformal class. As before, if $(r',s',v')$ are the components of a tangent vector, this vector is null if and only if the function
$$
r' \frac{\del f}{\del r} + s'\frac{\del f}{\del s} + v'
$$
has a zero.

In fact, the calculations are simpler if $f$ is not given exactly by the integral in Lemma \ref{formula for f}, but is instead adjusted by a certain constant which depends on $r$ and $s$. To describe this normalisation, fix a point $1$ in the original line $L$ and make a section for the $\C^2$-action on the space of lines by considering lines which meet $L$ at this fixed point. I.e., define the functions $f$ so that they vanish at $1$. 

On one of these lines, described as $\P^1_{r,s}$ as above, there are two marked points $0', \infty'$, where the $\C^2$-orbits are tangential, plus the point corresponding to $1$, where it meets $L$. So there is a preferred coordinate, $\tilde z$ say, on $\P^1_{r,s}$, normalised in the way indicated by the notation.

Now consider the solution of the jumping problem, regarded as depending on the parameters $r, s$. The derivative $\del f/\del r$ is initially a function on an open set in $E$. Since the jump data is independent of $r$ it actually extends to a meromorphic function on all of $E$ and, moreover, is pulled back from $\P^1_{r,s}$. Consider then $\del f/\del r$ and $\del f/\del s$ as meromorphic functions on $\P^1_{r,s}$. 

\begin{lemma}
There exist $M(r,s), N(r,s) \in \C^2$, such that 
$$
\frac{\del f}{\del r}
=
M(\tilde z ^{-1} - 1), 
\quad
\frac{\del f}{\del s}
=
N(\tilde z - 1).
$$
\end{lemma}
\begin{proof}
By choice of normalisation, the derivatives must vanish at $\tilde z =1$. Since the derivatives correspond to a section of $H^0(\mathcal O(1)\oplus \mathcal O(1))$ they can have only one zero and therefore also have a single pole which must be at either $0'$ or $\infty'$ since the orbits are transverse elsewhere. As $r$ varies with $s$ fixed, say, the geometry near $\infty'$ is unchanged; hence $\del f/ \del r$ is holomorphic near $\infty'$. Similarly $\del f/ \del s$ is holomorphic near $0'$.
\end{proof}

\begin{lemma}\label{metric formula}
Write $M_i, N_i$ for the components of $M, N$ and define two 1-forms on the space of lines by
$$
\rho_M = M_2 \diff v_1 - M_1 \diff v_2,
\qquad
\rho_N = N_2 \diff v_1 - N_1\diff v_2.
$$
Then the anti-self-dual conformal class on the space of lines is given by
$$
\diff r\, \rho_M - \diff s\,\rho_N
+
\frac{\rho_M \rho_N}{M_2N_1 - M_1N_2}
$$
\end{lemma}
\begin{proof}
The tangent vector $(r',s',v'_1,v'_2)$ is null if and only if the two quadratic equations
$$
r' M_i(1 - \tilde z) + s'N_i(\tilde z^2 - \tilde z) + v_i'\tilde z =0
$$
for $\tilde z$ have a simultaneous solution. Eliminating $\tilde z$ gives the result.
\end{proof}

Notice that the plane spanned by $\del/ \del r, \del/\del s$ is null, as it should be; the normalisation has been chosen so that this corresponds to the plane of lines passing through the central line at the fixed point $1$.
 
The problem then is to identify the functions $M$ and $N$, i.e., the residues of $\del f/\del r$ and $\del f/\del s$ when regarded as functions on $\P^1_{r,s}$. In principle, this is done by differentiating  $\int \phi(\tilde z) \alpha_q$ with respect to $r$ and $s$, where $\phi$ is regarded as a function on $\P^1_{r,s}$. The main difficulty here is that $\tilde z$ depends on $r, s$. 

Finding the coordinate $\tilde z$ amounts to uniformising $\P^1_{r,s}$. To clarify the $r, s$ dependence of $\tilde z$, suppose that $\P^1_{r,s}$ is uniformised in the sense that an identification $L \cong \P^1_{r,s}$ has been found. (Here, the central twistor line $L$ is playing the r\^ole of the ``standard'' Riemann sphere.) Let $w$ be a local coordinate near $0 \in L$  for which $\tau(w) = -w$; assume that the branch points corresponding to the parameter $r$ are $\pm \sqrt r$. The uniformising map is a biholomorphism $\tilde z = \tilde z_{r,s} \colon L \to \P^1_{r,s}$ such that near $0$, on the cut plane,
$$
\tilde z
=
\sum_{n=1}^\infty
c_n(w^2 - r)^{n/2},
$$
with $c_1 \neq 0$. Equally, the series can be inverted to give
$$
w^2 - r
=
\left(
\sum_{n=1}^\infty
d_n \tilde z^n
\right)^2,
$$
with $d_1 = c_1^{-1}$. 

Write $D(\tilde z) =\sum d_n \tilde z^n$. Of course, all the coefficients $c_n, d_n$ depend on the parameters $r, s$. Now, given the odd function $\phi(w)$ put $\phi(w) = wh(w^2)$ so that in terms of the $\tilde z$ coordinate,
$$
\phi (\tilde z) = \sqrt{ r + D(\tilde z)^2}H(\tilde z)
$$
where $H$ is the holomorphic function of $\tilde z$ (without cuts)
$$
H(\tilde z) = h(r + D(\tilde z)^2).
$$

Fix a point $z_0$ on $L$ and let $q$ be the corresponding point with coordinate $\tilde z_0 = \tilde z(z_0)$ on $\P^1_{r,s}$. Take the meromorphic 1-form $\alpha_q$ to be
$$
\alpha_q = \frac{\diff \tilde z}{\tilde z - \tilde z_0}.
$$
(That is, the fixed base point $q_*$ is the point $\tilde z = \infty$.) In terms of the coordinate $\tilde z$ then,  Lemma \ref{formula for f} gives
$$
f(q) = \frac{1}{4\pi i}
\int_\gamma 
\sqrt{r + D(\tilde z)^2} H(\tilde z) 
\frac{\diff \tilde z}{\tilde z - \tilde z_0} 
+ C(r,s).
$$
where $C$ is the normalising constant chosen so that $f(1)=0$.

Next, consider the derivative of $f(q)$, for fixed $q$, with respect to the parameter $r$. This appears very complicated, since in the formula above, $C$, $\tilde z_0$ and all the coefficients depend on $r$. However, thanks to the choice of normalisation, all that is needed is the residue of the derivative at $\tilde z_0 = 0$. In differentiating the formula above, the only term which contributes to the residue comes from differentiating $\tilde z_0$. That is, modulo terms which are bounded over $\tilde z_0 = 0$,
$$
\frac{\del f}{\del r}(q)
=
\frac{1}{4\pi i}
\frac{\del \tilde z_0}{\del r}
\int_\gamma
\sqrt{ r + D(\tilde z)^2} H(\tilde z) 
\frac{\diff \tilde z}{(\tilde z - \tilde z_0)^2}.
$$
Now the formula
$
\tilde z = c_1 \sqrt{w^2 -r} + \cdots
$ 
implies that 
\begin{eqnarray*}
\frac{\del \tilde z_0}{\del r}
&=&
\frac{-c_1}{2\sqrt{w_0^2 - r}} + \cdots\\
&=& 
\frac{-c_1^2}{2\tilde z_0} \quad \text{modulo bounded terms}
\end{eqnarray*}
Putting this together, the residue is $M= -c_1^2 I /8\pi i$ where $I$ is the contour integral
$$
I
=
\int_\gamma
\sqrt{r + D(\tilde z)^2} H(\tilde z) \frac{\diff \tilde z}{\tilde z^2}.
$$
An analogous argument gives a contour integral formula for $N$.

The essential difficulty in implementing this procedure to find the metric explicitly is the solution of the uniformisation problem for the Riemann surfaces $\P^1_{r,s}$, that is, finding the function $\tilde z$. While it follows from abstract theory that this exists, there is not an explicit procedure for finding it in general. That said, when more is known about the holomorphic data concerned, the uniformisation can often be found directly. This is the case for the following class of examples.

Let $d\geq2$ and $F_0$ be a rational function of degree $d$, thought of as a map from the Riemann sphere to itself. Suppose that $F_0(z) =z$ for $z=0,1,\infty$ and that $0$ and $\infty$ are critical values of $F_0$ with simple branching. Suppose that there are $2d-4$ other critical values $\zeta_1, \ldots,\zeta_{2d-4}$ in $\C^*$. The function $F_0$ gives an involution $\tau$ of small neighbourhoods of $0$ and  $\infty$ respectively, defined by interchanging the two local solutions of the equation $F_0(z)= w$. Moreover, suppose that $\lambda$ is a holomorphic function defined on neighbourhoods of $0$ and $\infty$ with odd leading term. Then the function $\phi(\tilde z) = \sqrt{\lambda (F_0(z))}$ is odd with respect to $\tau$ (the square root being well defined---up to an overall sign---because of the branching of $F_0$). In this way, the data $F_0, \lambda$ defines an anti-self-dual conformal class.

Now consider a smooth family of rational maps $F_{r,s}$ parametrised by $r$ near $0$ and $s$ near $\infty$ with the following properties:
\begin{itemize}
\item
$F_{0,\infty} = F_0$;
\item
$F_{r,s}(0)= r$, $F_{r,s}(1) = 1$,  $F_{r,s}(\infty) = s$;
\item
$F_{r,s}$ has critical values $r, s, \zeta_1, \ldots, \zeta_{2d-4}$
\end{itemize}
In other words, $F_{r,s}$ is obtained from $F_0$ by deforming the critical values at $0, \infty$ to $r, s$ respectively, whilst keeping the other critical values fixed. 

It follows from general theory that such a family exists. When $d=2, 3$ $F_{r,s}$ can be written down explicitly; for larger $d$, $F_{r,s}$ can be written down modulo the solutions of systems of algebraic equations.

The choice of critical values (and the fact that, as $F$ is a smooth family, the monodromy representations of $F_{r,s}$ and $F_0$ are the same) ensures that the composition $F_{r,s}^{-1} F_0$ can be defined as a holomorphic function on the plane with two small cuts near $0$ and $\infty$ and this precisely gives a solution of the uniformising problem for the Riemann surface $\P^1_{r,s}$ defined by the data. 

Following through the recipe above (for the derivative with respect to $r$) gives $d_1^2 =  F''_{r,s}(0)$, so that 
$$
M 
=
-\frac{1}{4 \pi i F''_{r,s}(0)}
\int_\gamma
\sqrt{ \lambda (F_{r,s}(\tilde z))}
\frac{\diff \tilde z}{\tilde z^2}.
$$
In particular, if $\lambda$ is a polynomial, the contour integrals involved have the form 
$$
\int_\gamma \sqrt{R(z)}\frac{\diff z}{z^2}
$$
where $R$ is a rational function. Thus we obtain solutions which can be written explicitly in  terms of classical special functions of this kind.

To illustrate this, consider the case when $d=2$. Here, the involution $\tau$ is globally defined, which is the situation considered previously in section \ref{when tau extends}. In the current notation, 
$$
F_{r,s}(z) = 
\frac{s(k z^2 +r)}{kz^2 +s}
$$
where $k=(1-r)/(1-s^{-1})$. The involution is $\tau(w) =w$ and $\phi(w) = \sqrt{\lambda(w^2)}$. The formula above gives
$$
M 
=
-\frac{s}{8\pi i k(s-r)}
\int_\gamma\sqrt{
\lambda\left(\frac{s(k z^2 +r)}{kz^2 +s}\right)}
\frac{\diff z}{z^2}
$$

Make the substitution 
$$
x = \sqrt{\frac{s(kz^2 + r)}{kz^2 + s}}
$$
Then \
$$
x \diff x 
=
z \diff z
\frac{ks (s - r)}{(kz^2 + s)^2},
$$
and
$$
kz^2 + s
=
\frac{s(r-s)}{x^2 - s}.
$$
This gives
$$
M(r,s)
= 
- \frac{1}{8 \pi i}
\sqrt{\frac{1-s}{1-r}}
\int_\gamma
\frac{x \phi(x)\, \diff x}
{(x^2 - r)^{3/2}(x^2-s)^{1/2}}.
$$
A similar calculation for the other coefficient gives
$$
N(r,s)
=
-\frac{1}{8 \pi i}
\sqrt{\frac{1-r}{1-s}}
\int_\gamma
\frac{x \phi(x)\, \diff x}
{(x^2 - r)^{1/2}(x^2-s)^{3/2}}.
$$

The fact that $M(r,s) = N(s,r)$ means that it is possible to change coordinates for the $\C^2$-action to make the metric in Lemma \ref{metric formula} explicitly surface othogonal. Doing  this recovers the formula of Theorem \ref{metric when tau extends}.

\bibliographystyle{alpha}
\bibliography{toric_asd}

\end{document}